\documentclass[12pt,twoside]{article}

\usepackage{amsmath,amsfonts,amsthm,amssymb,amsopn,amstext,amscd}
\usepackage{enumerate}
\usepackage{authblk}
\usepackage{graphicx}
\usepackage{float}
\usepackage{setspace} %For the space between lines
\usepackage{enumitem}
\usepackage{hyperref}
\allowdisplaybreaks

\newcommand{\cuad}{\begin{flushright}\vspace{-2ex}$\Box$\vspace{-2ex}\end{flushright}}

\newcommand{\fff}{{\mathcal F}}
\newcommand{\eee}{\mathcal{E}}

\cuad

\newtheorem{thm}{Theorem}
\newtheorem{definition}{Definition}
\newtheorem{proposition}[thm]{Proposition}

\newtheorem{lemma}[thm]{Lemma}

\begin{document}
\pagenumbering{arabic}
%\pagestyle{myheadings}
%\singlespace
%\onehalfspace

%%%%%%%%%%%%%%%%%%%%%%%%%%%%%%%%%%%%%%%%%%%%%%%%%%%%%%%%%%%%%%%%%%%%%
\title{An inhomogeneous controlled branching process\footnote{This is a plain preprint version of the following paper published in the journal \emph{Lithuanian Mathematical Journal}
(see the official journal website at \url{https://doi.org/10.1007/s10986-015-9265-0}):\\
González, M., Minuesa, C., Mota, M., del Puerto, I. \& Ramos, A. An inhomogeneous controlled branching process. \emph{Lithuanian Mathematical Journal} 55, 61--71 (2015). DOI: 10.1007/s10986-015-9265-0}}

\author[1,2]{Miguel Gonz\'alez}
\author[1,2]{Carmen Minuesa}
\author[1,2]{Manuel Mota}
\author[1,2]{In\'es del Puerto}
\author[1,3]{Alfonso Ramos}

\affil[1]{Department of Mathematics, Faculty of Sciences, University of Extremadura, Avda. Elvas, s/n, 06006, Badajoz, Spain}

\affil[2]{\url{mvelasco@unex.es}, \url{cminuesaa@unex.es}, \url{mota@unex.es}, \url{idelpuerto@unex.es}, \url{aramos@unex.es}}

\date{September 1, 2014}

\maketitle
%%%%%%%%%%%%%%%%%%%%%%%%%%%%%%%%%%%%%%%%%%%%%%%%%%%%%%%%%%%%%%%%%%%%%

 \begin{abstract}
A discrete time branching process where the offspring distribution is generation--dependent,
and the number of reproductive individuals is controlled by a random mechanism  is considered.
This model is a Markov chain but, in general, the
transition probabilities are non--stationary. Under not too
restrictive hypotheses, this model presents the classical duality
of branching processes: either  becomes extinct almost surely or
grows to infinity. Sufficient conditions for the almost
sure extinction and for a positive probability of indefinite
growth are provided. Finally  rates of growth of the process provided the non--extinction  are studied.
 \end{abstract}


%
%================= Text entry area ================================%
\section{Introduction}\label{s:1}

Branching processes are regarded as appropriate
probability models for the description of the extinction / growth of
populations (see \cite{haccou}).  The oldest and simplest discrete time branching process is the standard Bienaym\'e--Galton--Watson process that describes the evolution of a population in which each individual, independently of the others, gives rise to a random number of offspring
(in accordance with a common reproduction law), and then dies or is
not considered in the following counts. This standard model is not always adequate to describe actual phenomena, thus there are many variants of this  model  to deal with important properties of real--world populations. In particular, controlled
branching processes are useful to model some situations where some
kind of regulation is required.  Thus, for example, the existence of predators in the environment provokes that the population does not live in freedom, so that the survival of each animal (and therefore the possibility of
giving new births) will be strongly affected by this factor, and therefore it is required a control mechanism  at each generation that determines the number of progenitors in each generation that continues with the evolution of the population.

The development  of a controlled branching process consists of two phases: a reproductive
phase where individuals give birth to their offspring according to
a probability distribution, called reproduction law, and a control
phase in which is determined the number of potential progenitors of the generation. In this phase some individuals can be introduced into or removed from the population according to other probability distribution, called control law.

In the literature on controlled branching processes (see \cite{GonzalezMolinaPuertod} and  \cite{Yanev2}, and references therein), the control phase is
assumed to depend on the population size. On the other hand, in
the vast majority of works, the reproduction law is assumed to be
the same for every individual in any generation. However, it seems
reasonable to think that the reproductive abilities of the individuals of a
population may vary from one generation to another. One can find many published papers regarding standard or multitype Bienaym\'e--Galton--Watson processes whose reproduction laws vary with the generation, usually referred as varying environment models (see for example
\cite{Agresti},  \cite{Church}, \cite{SouzaBigginsa}  or
\cite{Fujimagari} for the standard  one and  \cite{BigginsCohnNerman},
\cite{HamblyJones} or \cite{Jones} for the  multitype one). But, until now, this possibility has not been considered in the class of the controlled branching processes, at least from a general viewpoint.

It is the aim of this paper to introduce and research the
controlled branching processes in varying environment. This model is defined as follows:

Let $\{X_{n,i}: n=0,1,\ldots; i=1,2,\ldots\}$ and $\{\phi_{n}(k):
n=0,1,\ldots; k=0,1,\dots\}$ be two independent sequences of
non--negative, integer--valued random variables satisfying:
\begin{itemize}
 \item[a)] The variables  $X_{n,i}$, $n=0,1,\ldots;
i=1,2,\ldots$, are independent and, for each $n$, $X_{n,i},
i=1,2,\dots$, have the same probability distribution,
$\{p_{n,j}\}_{j\ge0}$, with $p_{n,j}=P(X_{n,i}=j)$, $j\geq 0$, called {\it reproduction law} of the $n$th
generation.
 \item[b)] The stochastic processes $\{\phi_{n}(k)\}_{k\ge0}$, $n=0,1,\dots$ are assumed
to be independent and, for each $k$, the variables $\phi_{n}(k)$,
$n=0,1,\dots$,  have the same probability distribution, called the
{\it control law} for the population size $k$.
 \end{itemize}

The {\it controlled branching process in varying environment}
(CPVE) is a sequence of random variables, $\{Z_n\}_{n\ge0}$,
defined recursively by
\begin{equation}\label{edef}
Z_0=N,\quad Z_{n+1}=\sum_{i=1}^{\phi_n(Z_n)}X_{n,i}\ ,\quad
n=0,1,\ldots,
\end{equation}
 where the empty sum is defined to be 0, and $N$ an arbitrary non--negative integer.

Intuitively, $X_{n,i}$ represents the number of offspring produced
by the $i$th individual in the $n$th generation and $Z_n$
represents the total number of individuals in the $n$th
generation. Moreover, if $Z_n=k$, then $\phi_n(k)$ is the number
of progenitors in the $n$th generation that will produce their
offspring, according to the reproduction law $\{p_{n,j}\}_{j\geq
0}$. The offspring of these progenitors  forms the $(n+1)$st
generation of the population. The control is made in such a way that if
$\phi_n(k)>k$, new individuals are introduced into the population
and if $\phi_n(k)<k$, some individuals are removed from the
population.

The CPVE generalizes two
classical branching models widely studied in the scientific
literature on branching processes theory: if  $\phi_n(k)=k$,
$n=0,1,\dots, k=0,1,\dots$ we obtain the standard Bienaym\'e--Galton--Watson process in
varying environment. On the other hand, if the reproduction law is
the same for all the generations, i.e. $p_{n,j}$ only depends on
$j$, we obtain the controlled branching process with  random
control function (see for example \cite{GonzalezMolinaPuertod}, \cite{SebastyanovZubkov} and
\cite{Yanev} and references therein).

It is easy to prove that the CPVE is a Markov chain, in general inhomogeneous. Our objectives in this paper are to establish the basic properties of this model and to study its long--term behaviour. For that, besides this introduction, the paper is organized as follows. In
Section \ref{s:2}, conditions for the  extinction--explosion duality to hold are stated and the extinction problem is tackled.
Section \ref{s:3} is devoted to studying the rate of convergence of the process on the non--extinction set.
The proofs  are relegated to Section
\ref{s:4}, in order to make easier  the reading of the paper.

\section{The extinction problem }\label{s:2}

Homogenous branching processes  often show a dual long term
behaviour: either become extinct or grow to infinity. However, to
obtain this behaviour in  inhomogenous processes, additional
regularity conditions are required.  In the following result we
provide sufficient conditions, given in terms of the reproduction
and control laws, for the CPVE to present also this duality.

\begin{thm}\label{thm1}
Let $\{Z_n\}_{n\ge0}$ be a CPVE satisfying:
\begin{enumerate}
    \item [(i)] $P(\phi_0(0)=0)=1$.
    \item [(ii)] $\liminf_{n\to\infty}p_{n,0}>0$.
\end{enumerate}
Then
\begin{equation}\label{ecc1}
P(Z_n\to 0)+P(Z_n\to\infty)=1.
\end{equation}

\end{thm}

Condition (i) in Theorem \ref{thm1} means that $0$ is an absorbing state. Whereas condition (ii) is trivially verified  for many families of offspring
distributions, for example, geometric probability distributions with parameter $\alpha_n$,
being $0<\liminf_{n\to\infty}\alpha_n<1$.

At this point, it is
worth mentioning  that for the Galton--Watson process in varying environment $\sum_{n=0}^\infty(1-p_{n,1})=\infty$ is
a sufficient condition for the duality extinction--explosion to
hold (see \cite{Lindvall}). Condition (ii) of the previous
theorem is stronger. In fact, if $\liminf_{n\to\infty}p_{n,0}>0$
then $\liminf_{n\to\infty}(1-p_{n,1})>0$ and therefore
$\sum_{n=0}^\infty(1-p_{n,1})=\infty$. The presence of the control variables
makes difficult to provide  a sharper condition than (ii)  or a necessary and sufficient condition for (\ref{ecc1}) to hold.

We are now interested in stating sufficient conditions for the
almost sure extinction and for the indefinite growth of the CPVE.
Such conditions will be given in terms of the first and second order
moments of the reproduction and control laws and in terms of their respective
probability generating functions. Let us introduce the notation. For $n,k=0,1,\dots$, define:
$$
m_n=E[X_{n,1}]\ ,\quad \sigma_n^2=\mbox{Var}[X_{n,1}]\ ,\quad
\eee(k)=E[\phi_0(k)]\ ,\quad \tau^2(k)=\mbox{Var}[\phi_0(k)]\
$$(assumed finite) and
$$
f_n(s)=E[s^{X_{n,1}}],\quad g_{k}(s)=E[s^{\phi_0(k)}],\quad
\ 0\le s\le 1.
$$
From (\ref{edef}) it follows that, for $n,k=0,1,\dots$,
\begin{eqnarray}
E[Z_{n+1}|Z_n=k]&=&m_n\eee(k),\label{e1}\\
\mbox{Var}[Z_{n+1}|Z_n=k]&=&m_n^2\tau^2(k)+\sigma_n^2\eee(k).\label{e2}
\end{eqnarray}

As usual, let us denote the probability of extinction as $q=P(Z_n \to 0)$.
We will assume throughout the paper that $P(\phi_0(0)=0)=1$. Under this
hypothesis the extinction probability can be rewritten as
$q=\lim_{n\to\infty}P(Z_n=0)$. We will   also assume
$\liminf_{n\to\infty}p_{n,0}>0$, so that Theorem \ref{thm1} applies and (\ref{ecc1}) holds.
Consequently, we deduce that
$P(Z_n\to\infty)=1-q=\lim_{n\to\infty}P(Z_n>0)$. First, in the
next two results, we will provide some sufficient conditions for
the almost sure extinction of the process.
\begin{thm}\label{thm2}
Let $\{Z_n\}_{n\ge0}$ be a  CPVE such that
\begin{equation}\label{ecn4}
\limsup_{k\to\infty}k^{-1}\eee(k)<\liminf_{n\to\infty}m_n^{-1}.
\end{equation}

Then $q=1$.
\end{thm}
\begin{thm}\label{thm3}
Let $\{Z_n\}_{n\ge0}$ be a CPVE and let
$\gamma_n(s)=\inf_{k\ge1}g_k(f_n(s))$, $n=0,1,\ldots$, $0\leq
s<1$. If for some $s$, $0\leq s<1$,\begin{equation}\label{new}
\liminf_{n\to\infty}\gamma_n(s)>0,
\end{equation}
then $q=1$.
\end{thm}

Notice that Theorem \ref{thm2} includes as particular cases the
conditions given in \cite{Fujimagari} for the almost sure
extinction of the standard Bienaym\'e--Galton--Watson branching process in varying
environment  and those given in \cite{GonzalezMolinaPuertoa} for the almost sure
extinction of  the controlled
branching processes with random control function.

\vspace*{0.5cm}In order to provide sufficient conditions for a
positive probability of non--extinction, we provide the definition
of uniformly supercritical CPVE.

\begin{definition}
A CPVE is said to be {\em uniformly supercritical} if there exists
a constant $\eta>1$ such that
\begin{equation}\label{eq2.5}
\liminf_{n\to\infty}m_n\ge\eta\,k\eee(k)^{-1}\quad\mbox{ for
every $k\ge1$}.
\end{equation}
\end{definition}

In  $\eee(k)=k$,  for every $k$, then   (\ref{eq2.5}) is a sufficient condition
for a Galton--Watson branching process in
varying environment to be uniformly supercritical according to the definition given in  \cite{SouzaBigginsa}.

An uniformly supercritical CPVE has a positive probability of
indefinite growth if it satisfies some conditions on the second
order moments of the offspring and control laws. We can establish
two results. The first one uses a similar methodology to that applied  in
\cite{SouzaBigginsa} for the Galton--Watson branching process in
varying environment and the second  generalizes the conditions given  in
\cite{GonzalezMolinaPuertoa} for  the controlled branching
process with random control function.

\begin{thm}\label{thm4}
Let $\{Z_n\}_{n\ge0}$ be an uniformly supercritical CPVE and let
$\eta>1$ satisfying (\ref{eq2.5}). Assume also that the following
conditions hold:
\begin{enumerate}
    \item [(i)] There exists a constant $\gamma>0$ such that
    $\gamma/\eta^2<1$ and
    $$
    \limsup_{n\to\infty}m_n^2\le \gamma\,\frac{ k^2}{d^2(k)}\quad\mbox{ for every } k\ge 1,
    $$
    with $d^2(k)=E[\phi_0^2(k)]$, $k=0,1,\dots$
    \item [(ii)]
    $$\sum_{n=0}^\infty\frac{\sigma_n^2}{m_n^2\eta^n}<\infty\,.$$
\end{enumerate}
Then $q<1$.
\end{thm}

\begin{thm}\label{thm5}
Let $\{Z_n\}_{n\ge0}$ be an uniformly supercritical CPVE and let
$\eta>1$ satisfying (\ref{eq2.5}). Assume also that the sequences
$\{\eee(k)/k\}_{k\ge1}$ and $\{\tau^2(k)/k\}_{k\ge1}$ are bounded
and that there exists $\delta>0$ such that
$$
 \sum_{n=0}^\infty\frac{m_n^2+\sigma^2_n}{(\eta-\delta)^n}<\infty.
$$
Then $q<1$.
\end{thm}
From Theorems \ref{thm2}, \ref{thm4} and \ref{thm5}, it is deduced that it seems quite important the behaviour of the sequence of the expected growth rates per individual when, in  a certain generation, there are $k$ individuals, that is,  $E[Z_{n+1}Z_n^{-1}\mid
Z_n=k]=m_nk^{-1}\eee(k)$, $k=1,2,\ldots$, in order to determine the extinction probability. Indeed, this fact is usual in most of the branching models. Again, it is not surprising that its behaviour respect to the value 1  establishes somehow the threshold  for extinction or non--extinction of the process. In fact, conditions (\ref{ecn4}) and (\ref{eq2.5}) can be rewritten  as $m_nk^{-1}\eee(k)<1$, for all $k\geq k_0$ and $n\geq n_0$, $k_0,n_0>0$ and $m_nk^{-1}\eee(k)>1$, for all $n\geq N_0$, for all $k\geq 1$, $N_{0}>0$, respectively. It is a matter for further research to study the behaviour of the process when this double indexed sequence $\{m_nk^{-1}\eee(k)\}_{n, k\geq 0}$ approaches to 1.

\section{Asymptotic behaviour}
\label{s:3}
If $q<1$, we are interested in the rate of  growth of $\{Z_n\}_{n\geq 0}$ on the non--extinction set. In particular, we wonder if there exist sequences of positive constants $\{r_n\}_{n\geq 0}$ such that $\lim_{n\to\infty} Z_n/r_n$ exists almost surely and $P(0<\lim_{n\to\infty} Z_n/r_n<\infty)>0$.
To this end, let us  assume the asymptotic linear growth of the mathematical expectations of the control means, i.e.,  suppose that $\tau=\lim_{k \to \infty} k^{-1}\eee(k)$ exists and is finite. Let us consider the sequences $$r_n=\tau^n\prod_{i=0}^{n-1}m_i\quad \mbox{ and }\quad W_n=r_n^{-1}Z_n.$$

By the supermartingale convergence theorem, if
$\{k^{-1}\eee(k)\}_{k \geq 1}$ is a monotonic increasing
sequence, then $\{W_n\}_{n\geq 0}$ converges almost surely  to a non--negative and finite random variable, $W$, as $n\to\infty$. So, from now on, we assume that the sequence $\{k^{-1}\eee(k)\}_{k \geq 1}$ is  increasing and let denote $\delta_k=\tau-k^{-1}\eee(k)$. Consequently $\{\delta_k\}_{k\geq 0}$ is a non--increasing sequence with limit equal to zero. Moreover we assume that the process $\{Z_n\}_{n\geq 0}$ is uniformly supercritical so that $\tau\liminf_{n\to\infty}m_n\geq \eta $, for some constant $\eta>1$, i.e., there exists $n_0$ such that for all $n\geq n_0$, $\tau m_n\geq \eta $. For simplicity, we will assume without loss of generality  that $n_0=0$.

The following result establishes a condition for the existence of the
limit of $\{E[W_n]\}_{n\geq 0}$ as $n\to \infty$. Such a limit
will be positive and finite if the process starts with a large
enough number of individuals.
\begin{proposition}\label{prop1}
Let $\{Z_n\}_{n\ge0}$ be an uniformly supercritical CPVE and  let
$\eta>1$ satisfying (\ref{eq2.5}). Assume that the sequence $\{\delta_k\}_{k\ge1}$
is a non--increasing sequence and  $\sum_{k=1}^\infty k^{-1}\delta_k<\infty$.
Then there exists $N_0$ such that $$\lim_{n\to\infty}E[W_n]>0,\ \mbox{ if }Z_0=N>N_0 \mbox{ with } q<1.$$
\end{proposition}

Finally we provide the following result where we prove convergence almost sure, in $L^1$ and in $L^2$ of $\{W_n\}_{n\geq 0}$ to a non--degenerate  random variable by assuming that $q<1$. Therefore we establish the geometric growth of the CPVE in the  uniformly supercritical cases.
\begin{thm}\label{thm6}
Let $\{Z_n\}_{n\ge0}$ be an uniformly supercritical CPVE and let
$\eta>1$ satisfying (\ref{eq2.5}). Assume that:
\begin{itemize}
  \item[(i)] The sequence $\{\delta_k\}_{k\ge1}$ is non--increasing and $\sum_{k=1}^\infty k^{-1}\delta_k<\infty$.
  \item[(ii)] The sequence $\{k^{-2}\tau^2(k)\}_{k\ge1}$ is non--increasing and $\sum_{k=1}^\infty k^{-3}\tau^2(k)<\infty$.
  \item[(iii)] $\displaystyle\sum_{n=0}^\infty\nolimits\frac{\sigma_n^2}{m_n^2\eta^n}<\infty$.
\end{itemize}
Then $\{W_n\}_{n\ge0}$ is a $L^2$--bounded supermartingale and converges almost surely, in $L^1$ and in $L^2$ to the random variable $W$ that is finite almost surely and non--degenerate at 0.
\end{thm}

\section*{Concluding remark}
We have introduced a new branching model called the CPVE presenting the novelty of  joining the possibility that the reproduction laws vary with the generation  to the incorporation of a random mechanism that determines the number of progenitors in each generation. The methodologies developed independently for the controlled branching processes and for the Galton-Watson processes in varying environment have been adapted to study the extinction problem and the rates of growth for the CPVE. It is interesting to mention that a CPVE could be also thought as a general branching process with size and time-dependent reproduction laws in a non-trivial way as follows: Using the notation in (\ref{edef}) $$Z_{n+1}=\sum_{i=1}^{Z_n}Y_{n,i}(Z_n), n=0,1,\ldots,\quad \mbox{ with }$$
$$Y_{n,i}(Z_n)=X_{n,i}+Z_n^{-1}\sum_{j=1}^{\phi_n(Z_n)-Z_n}X_{n,j}I_{\{\phi_n(Z_n)>Z_n\}}-
Z_n^{-1}\sum_{j=1}^{Z_n-\phi_n(Z_n)}X_{n,j}I_{\{\phi_n(Z_n)<Z_n\}}.$$
In this case, the expected value  $\mu_{n, k}=E[Y_{n,i}(Z_n)\mid Z_n=k]=k^{-1}m_n\eee(k)$, $n,k=1,2,\ldots$, depends on the generation  and the population size.
Therefore, one could expect that adapting techniques of size-dependent branching processes, it would be possible to obtain some results about this process (indeed Theorem \ref{thm6} follows these ideas). This is an interesting open topic for further research. However this approach implies that the mathematical modelling of the control achieved on the population sizes  at each generation could be diluted. In this sense, as is proposed in the present paper, it is worth putting together the ideas of control on the population and generation dependent reproduction law in such a way that both features appear explicitly in the definition of the model. This allows us to research regularity conditions for the control and reproduction laws that lead to the extinction or survival of the process. It also provides us a greater capacity of observation to face up to the estimation problem of the parameters of the model and therefore to develop more easily some relevant potential applications of the process.

\section{Proofs}\label{s:4}
\emph{Proof of Theorem \ref{thm1}}
Since (ii) holds, there exist $n_1>0$ and $0<a<1$ such that
$p_{n,0}\ge a >0$ for every $n\ge n_1$. Moreover, $g_k(a)>0$ for
every $k=0,1,\dots$

First, given $i, j>0$, let us prove that for every $m>0$
\begin{eqnarray}\label{eq1}
&&\hspace*{-2.75cm}\nonumber P(Z_{n+k}=j\ \mbox{ for at least $m$
values of }\ k>0|Z_n=i)
\\&\le&(1-g_j(a))^{m-1}(1-g_i(a))\quad\mbox{for all $n\ge n_1$.}
\end{eqnarray}
We proceed by induction on $m$.

For $m=1$, taking into account that, by virtue of (i), $Z_{n+1}=0$
implies that $Z_{n+k}=0$  for all $k\ge1$, we have
that, for all $n\ge n_1$,
\begin{eqnarray*}
&&\hspace*{-0.75cm}P(Z_{n+k}=j\ \mbox{ for some }\ k>0|Z_n=i)\le
1-P(Z_{n+1}=0|Z_n=i)\\&&=1-P\left(\sum_{l=1}^{\phi_n(i)}X_{n,l}=0\right)
=1-E[p_{n,0}^{\phi_n(i)}]=1-g_i(p_{n,0})\le 1-g_i(a).
\end{eqnarray*}

Suppose that (\ref{eq1}) holds for all the positive integers less
than or equal to $m$ and let us prove it for $m+1$. Then, taking
again $n\ge n_1$ and applying Markov property, we obtain
\begin{eqnarray*}
&&\hspace*{-1.75cm}P(Z_{n+k}=j\ \mbox{ for at least $m+1$ values of }\ k>0|Z_n=i)\\
%=\,&\sum_{k=1}^\infty P(Z_{n+l}\ne j,\ 0<l<k,\ Z_{n+k}=j,\
%Z_{n+h}=j\ \mbox{ for at least $m$ values of }\ h>k|Z_n=i)\\
&=\,&\sum_{k=1}^\infty P(Z_{n+h}=j\ \mbox{ for at least $m$ values
of }\ h>k|Z_{n+k}=j)\\&&\phantom{\sum_{k=1}^\infty}\cdot
P(Z_{n+k}=j,\ Z_{n+l}\ne j,\ 0<l<k|Z_n=i)\\
&\le\,&(1-g_j(a))^m \sum_{k=1}^\infty P(Z_{n+k}=j,\ Z_{n+l}\ne j,\
0<l<k|Z_n=i)\\
&=\,&(1-g_j(a))^mP(Z_{n+k}=j\ \mbox{ for some }\ k>0|Z_n=i)\\
&\le&(1-g_j(a))^m(1-g_i(a)).
\end{eqnarray*}

Hence, since $g_j(a)>0$ for every $j>0$, we have
 $$
P(Z_n=j\quad \mbox{i.o.})=\sum_{i=1}^\infty P(Z_n=j\quad \mbox{
i.o. }, n\ge n_1 |Z_{n_1}=i)P(Z_{n_1}=i).
 $$
But
 \begin{eqnarray*}
&& \hspace*{-1.75cm}P(Z_n=j\quad \mbox{ i.o. }, n\ge n_1
|Z_{n_1}=i)
 \\&=&
\lim_{m\to\infty}P(Z_{n_1+k}=j\ \mbox{ for at least $m$ values of
}\ k>0|Z_{n_1}=i)  \\&\le& \lim_{m\to\infty}
(1-g_j(a))^m(1-g_i(a))=0.
 \end{eqnarray*}
 So $P(Z_n=j\quad \mbox{i.o.})=0$.

Finally, taking into account that there exists an integer $B>0$
such that $(\{Z_n\to0\}\cup\{Z_n\to\infty\})^c\subseteq \{0<Z_n\le
B\quad \mbox{i.o. }\}$ and that
$$
P(0<Z_n\le B\quad \mbox{i.o. })=\sum_{j=1}^BP(Z_n=j\quad \mbox{i.o.})=0,
$$
we conclude the proof.
\vspace*{0.25cm}

\emph{Proof of Theorem \ref{thm2}}
We need the following auxiliary lemma, for which we provide a
sketch of the proof, the details can be found in \cite{Martinez},
p.\ 41:

\begin{lemma}\label{lemma1}
Let $\{X_n\}_{n\ge0}$ be a sequence of non negative random
variables and $\{\fff_n\}_{n\ge0}$ a sequence of $\sigma$--algebras
such that $X_n$ is $\fff_n$--measurable for all $n$. If there
exists a constant $A>0$ such that, for every $n$,
$E[X_{n+1}|\fff_n]\le X_n$ almost surely on $\{X_n\ge A\}$, then
$P(X_n\to\infty)=0$.
\end{lemma}

\noindent {\bf Proof of the Lemma:}

Let $A>0$ satisfying the hypothesis of the lemma. It is enough to
prove that, for every $N>0$, $ P\left(\inf_{n\geq N} X_n\geq
A,X_n\to\infty\right)=0. $ Fixed $N>0$, define the {\it stoping
time} $T(A)$ by $\inf\{n\geq N: X_n< A \}$  if $\inf_{n\geq
N}X_n< A$ and by $\infty$ otherwise. Define also the sequence of
random variables $\{Y_n\}_{n\geq 0}$, with $Y_n$ for $n\geq 0$ as
follows
$$
Y_n=\left\{\begin{array}{cc} X_{N+n} & \mbox{ if } N+n\leq T(A), \\
X_{T(A)} & \mbox{ if } N+n> T(A). \end{array}\right.
$$

Since $E[X_{n+1}|\fff_n]\le X_n$  almost surely on $\{X_n\ge A\}$,
$\{Y_n\}_{n\geq 0}$ is a non--negative supermartingale
and applying the martingale convergence theorem, we obtain the
almost sure convergence of the sequence $\{Y_n\}_{n\geq 0}$ to a
non--negative and finite limit, and therefore the proof of the
Lemma ends.

\vspace{.4cm}

Let us prove the theorem. By hypothesis there exist $A>0$ and
$n_0>0$ such that
\begin{equation}\label{eq2}
k^{-1}\eee(k)<m_n^{-1},\qquad \mbox{ for all } k\ge A\ \mbox{and}\
n\ge n_0.
\end{equation}
Assume without loss of generality that $n_0=0$. Otherwise we would
proceed with the sequence $\{Z_n\}_{n\ge n_0}$.

Denote by $\fff_n=\sigma(Z_0,\dots,Z_n)$, $n=0,1,\dots,$ i.e. the $\sigma$--algebra generated by the random variables $\{Z_0,\ldots,Z_n\}$. Since
$\{Z_n\}_{n\ge0}$ is a Markov chain, and using
(\ref{e1}) and (\ref{eq2}), we deduce that, for all $n$,
$$
E[Z_{n+1}|\fff_n]=E[Z_{n+1}|Z_n]=m_n\eee(Z_n)\le Z_n\quad\mbox{ on
}\ \{Z_n\ge A\}.
$$
Now, since we are assuming that $P(Z_n\to\infty)=1-q$, applying
Lemma \ref{lemma1} the proof is finished.

\vspace*{0.25cm}

\emph{Proof of Theorem \ref{thm3}}
Let us denote the probability generating function of $Z_n$ by $F_n(s)$, $0\leq s\leq 1$. We have, for $n=0,1,\ldots$, and
$0\leq s<1$,
$$
F_{n+1}(s)=E[g_{Z_{n}}(f_n(s))] = P(Z_n=0) +
E[g_{Z_{n}}(f_n(s))I_{\{Z_n>0\}}].
$$
Using (\ref{ecc1}),
$\lim_{n\to\infty}s^{Z_n}=I_{\{Z_n\to 0\}}$ almost surely and
therefore $\lim_{n\to\infty}F_{n+1}(s)=q$.

Hence, for all $0\leq s<1$,
$$
\lim_{n\to\infty}E[g_{Z_{n}}(f_n(s))I_{\{Z_n>0\}}]=0,
$$
and by Fatou's lemma,
$$
E[\liminf_{n\to\infty}g_{Z_{n}}(f_n(s))I_{\{Z_n>0\}}]=0.
$$

Now on $\{Z_n\to\infty\}$,
$$
\liminf_{n\to\infty}g_{Z_{n}}(f_n(s))I_{\{Z_n>0\}}\geq
\liminf_{n\to\infty}\gamma_n(s),
$$
and therefore $0=(1-q)\liminf_{n\to\infty}\gamma_n(s)$ for all
$0\leq s<1$, so using (\ref{new}) it is deduced $q=1$.
\vspace*{0.25cm}

\emph{Proof of Theorem \ref{thm4}}
Since the CPVE is uniformly supercritical and hypothesis (i) holds, we can
take $n_0$ such that for all $n\ge n_0$
\begin{equation}\label{eq3}
\eee(k)\ge m_n^{-1}\eta k\quad\mbox{  for every $k\ge 1$}
\end{equation}
and
\begin{equation}\label{eq3.5}
d^2(k)\le m_n^{-2}\gamma k^2 \quad\mbox{  for every $k\ge 1$}.
\end{equation}
Assume without loss of generality that $n_0=0$. Otherwise we would
proceed with the sequence $\{Z_n\}_{n\ge n_0}$.

We will make use of the fact that, for every non--negative random
variable $Y$, the following inequality holds:
\begin{equation}\label{eq4}
P(Y>0)\ge\left(\frac{\mbox{Var}[Y]}{E[Y]^2}+1\right)^{-1}.
\end{equation}
 Let us denote $T_n=\phi_n(Z_n)$. Using  (\ref{e1}), (\ref{e2}) and $$\mbox{Var}[Z_{n+1}]=\mbox{Var}[E[Z_{n+1}|Z_n]]+E[\mbox{Var}[Z_{n+1}|Z_n]],$$ we obtain
$$
\frac{\mbox{Var}[Z_{n+1}]}{E[Z_{n+1}]^2}=\frac{\mbox{Var}[T_n]}{E[T_n]^2}+\frac{\sigma_n^2}{m_n^2E[T_n]}
=\frac{E[T^2_n]}{E[T_n]^2}-1+\frac{\sigma_n^2}{m_n^2E[T_n]}\,.
$$
Using recursively (\ref{eq3}),  we
have
$$
E[T_n]=E[\eee(Z_n)]\ge m_n^{-1}\eta E[Z_n]= \eta E[T_{n-1}]\ge
E[T_0]\eta^n.
$$
Moreover, using (\ref{eq3.5}),
$$
E[T_n^2]=E[d^2(Z_n)]\le m_n^{-2}\gamma E[Z_n^2].
$$
Since $\gamma/\eta^2<1$, we deduce from the previous equations
that
$$
\frac{\mbox{Var}[Z_{n+1}]}{E[Z_{n+1}]^2}\le
\frac{E[Z^2_n]\gamma}{E[Z_n]^2\eta^2}-1+\frac{\sigma_n^2}{m_n^2E[T_0]\eta^n}\le
\frac{\mbox{Var}[Z_n]}{E[Z_n]^2}+\frac{\sigma_n^2}{m_n^2E[T_0]\eta^n}.
$$
By iteration we obtain that, for every $n\ge0$,
$$
\frac{\mbox{Var}[Z_{n+1}]}{E[Z_{n+1}]^2}\le
\frac{1}{E[T_0]}\sum_{j=0}^\infty\frac{\sigma_j^2}{m_j^2\eta^j},
$$
and, by hypothesis (ii), this series is convergent. So, applying (\ref{eq4}),
for every $n\ge0$ we obtain
$$
P(Z_{n+1}>0)\ge\left(\frac{1}{E[T_0]}\sum_{j=0}^\infty\frac{\sigma_j^2}{m_j^2\eta^j}+1\right)^{-1}.
$$
The right hand side of this inequality is positive and does not
depend on $n$, therefore
$$
P(Z_n\to\infty)=\lim_{n\to\infty}P(Z_n>0)
\ge\left(\frac{1}{E[T_0]}\sum_{j=0}^\infty\frac{\sigma_j^2}{m_j^2\eta^j}+1\right)^{-1}>0,
$$
which finishes the proof.
\vspace*{0.25cm}

\emph{Proof of Theorem \ref{thm5}}
We will prove that $P(Z_n\to\infty)=1-q>0$. Since the process is
uniformly supercritical, there exists $n_0$ such that for all
$n\ge n_0$ and for every $k\ge 1$
\begin{equation}\label{eq5}
m_n\,\frac{\eee(k)}{k}\ge \eta.
\end{equation}
Assume without loss of generality that $n_0=0$. Otherwise we would
proceed with the sequence $\{Z_n\}_{n\ge n_0}$, by showing that
$P(Z_n\to\infty|Z_{n_0}=N')>0$ for some $N'>0$.

Take $\delta'>0$ such that $\delta'<\min\{\eta-1,\delta\}$ and
denote $A_n=\{Z_{n+1}>(\eta-\delta')Z_n\}$. Since
$\eta-\delta'>1$, it is immediate that $\cap_{n=0}^\infty
A_n\subseteq \{Z_n\to\infty\}$, so it is enough to prove that
$P(\cap_{n=0}^\infty A_n)>0$.

By the Markov property and using that $Z_0=N$ we have
\begin{eqnarray}\label{eq6}\nonumber
P\left(\bigcap_{n=0}^\infty
A_n\right)=&&\lim\limits_{l\to\infty}P\left(\bigcap_{n=0}^l
A_n\right)
=P(A_0)\lim\limits_{l\to\infty}\prod_{n=1}^lP\left(A_n\left|\right.\bigcap_{j=0}^{n-1}
A_j\right)\\\ge\,&&
P(A_0)\prod_{n=1}^\infty\inf\limits_{k>(\eta-\delta')^nN}P(A_n|Z_n=k).
\end{eqnarray}
Take $a$ and $b$, bounds for the sequences $\{\eee(k)/k\}_{k\ge1}$
and $\{\tau^2(k)/k\}_{k\ge1}$, respectively. Applying (\ref{e2}), (\ref{eq5}) and Chebyshev's inequality we obtain
\begin{eqnarray*}
P(A_n^c|Z_n=k)&\le& P(Z_{n+1}\le m_n\eee(k)-k\delta'|Z_n=k)\\&\le&
P(|Z_{n+1}-m_n\eee(k)|\ge k\delta'\mid Z_n=k)\\
&\le\,&\frac{\mbox{Var}[Z_{n+1}|Z_n=k]}{k^2\delta'^2}=\frac{m_n^2\tau^2(k)+\sigma_n^2\eee(k)}{k^2\delta'^2}
\le \frac{m_n^2 a+\sigma_n^2 b}{k\delta'^2}.
\end{eqnarray*}
Therefore, from (\ref{eq6}),
\begin{eqnarray*}
P\left(\bigcap_{n=0}^\infty A_n\right)\ge\,&
P(A_0)\prod\limits_{n=1}^\infty\inf\limits_{k>(\eta-\delta')^nN}\left(1-
\displaystyle\frac{m_n^2 a+\sigma_n^2 b}{k\delta'^2}\right)\\
=\,&P(A_0)\prod\limits_{n=1}^\infty\left(1-
\displaystyle\frac{m_n^2 a+\sigma_n^2
b}{(\eta-\delta')^nN\delta'^2}\right).
\end{eqnarray*}
Since $\delta'<\delta$, by hypothesis the series
$\sum_{n=1}^\infty m_n^2/(\eta-\delta')^n$ and $\sum_{n=1}^\infty
\sigma_n^2/(\eta-\delta')^n$ are convergent and consequently
$$
P\left(\bigcap_{n=0}^\infty A_n\right)\ge
P(A_0)\prod_{n=1}^\infty\left(1- \frac{m_n^2 a+\sigma_n^2
b}{(\eta-\delta')^nN\delta'^2}\right)>0,
$$
which finishes the proof.
\vspace*{0.25cm}

\emph{Proof of Proposition \ref{prop1}}
In \cite{Kl85}, it was proved that under the hypotheses satisfied by the  sequence $\{\delta_k\}_{k\geq 0}$, there exists a positive and non--increasing  function $\delta(x)$ such that $\delta_k\leq \delta(k)$, for all $k$, $x\delta(x)$ is concave and  $\sum_{k=1}^\infty k^{-1}\delta(k)<\infty$. Thus, by Jensen's inequality one can check that $$0\leq E[W_n]-E[W_{n+1}]\leq \tau^{-1}E[W_n]\delta(\eta^nE[W_n]),$$.

Now, one can use Lemma 2 in \cite{Kl84}, by assuming $a_n=E[W_n]$, $f=\tau^{-1}\delta$ and $m=\eta$, to conclude the result.
\vspace*{0.25cm}

\emph{Proof of Theorem \ref{thm6}}
Under the hypotheses of the theorem one has that $\{W_n\}_{n\geq 0}$ is a supermartingale.
It will be enough to check that it is $L^2$--bounded to obtain its $L^1$--convergence to $W$.
Moreover the limit $W$ is non--degenerate at 0 because  under $L^1$--convergence $\lim_{n\to\infty}E[W_n]=E[W]$, and using  Proposition \ref{prop1}, this limit is greater than 0 if $N$ is large enough.

Let us prove that $\{E[W_n^2]\}_{n\geq 0} $ is a bounded sequence.
The proof  is a suitably adaptation of Theorem 3 in \cite{GonzalezMolinaPuertoc}, so that we only show the main steps.
Some calculations lead us to $$E[W_{n+1}^2]= E[W_n^2]+\frac{1}{\tau
^2}E\left[W_n^2\left(\frac{\tau^2(Z_n)}{Z_n^2}+\delta_{Z_n}^2-2\tau
\delta_{Z_n}\right)\right ]+\sigma_n^2\frac{E[\eee(Z_n)]}{r
^2_{n+1}}.$$

By considering the properties of the sequences $\{\delta_k\}_{k\geq 0}$ and $\{\tau^2(k)\}_{k\geq 0}$ and
results in \cite{GonzalezMolinaPuertoc}, there exist positive and non--increasing functions $\delta(x)$ and $h(x)$
such that:

\begin{itemize}
\item[a)] $\delta_k\leq \delta(k)$, for all $k\geq 0$, $\sum_{k=1}^\infty k^{-1}\delta(k)<\infty$ and the functions $x\delta(x)$, $x\delta(x^{1/2})$ and $x\delta^2(x^{1/2})$ are concave.
\item[b)] $k^{-2}\tau^2(k)\leq h(k)$, for all $k\geq 0$, $\sum_{k=1}^\infty k^{-1}h(k)<\infty$ and the function $xh(x^{1/2})$ is concave.
\end{itemize}

Therefore,
$$E[W_{n+1}^2]\leq E[W_n^2]\left (1+\frac{1}{\tau^2}(h(E[Z_n])+\delta^2(E[Z_n])+2\tau
\delta(E[Z_n]))\right)+\sigma_n^2\frac{E[\eee(Z_n)]}{r
^2_{n+1}}.$$

By Proposition \ref{prop1}, there exists $c>0$ such that $E[W_n]>c$ for all $n$, and since
the process is uniformly supercritical, one deduces that $E[Z_n]\geq c\eta^n$.
Hence, using that $h(x),\ \delta(x)$ and $\delta^2(x)$ are non--increasing functions, one obtains that

\begin{eqnarray*}E[W_{n+1}^2]&\leq& Z_0^2\prod_{i=0}^n  \left(1+\frac{1}{\tau^2}(h(c\eta^i)+\delta^2(c\eta^i)+2\tau
\delta(c\eta^i))\right)\\ &+& \sum_{i=0}^n\sigma_i^2\frac{E[\eee(Z_i)]}{r
^2_{i+1}}\prod_{j=i+1}^n  \left(1+\frac{1}{\tau^2}(h(c\eta^j)+\delta^2(c\eta^j)+2\tau
\delta(c\eta^j))\right).
\end{eqnarray*}

To conclude that $\{E[W_n^2]\}_{n\geq 0}$ is a bounded sequence, it is enough to check that
\begin{equation}\label{ecn3}
\sum_{n=0}^\infty \sigma_n^2\frac{E[\eee(Z_n)]}{r
^2_{n+1}}< \infty
\end{equation}
and

\begin{equation}\label{ecn2}
\prod_{n=0}^\infty  \left(1+\frac{1}{\tau^2}(h(c\eta^n)+\delta^2(c\eta^n)+2\tau
\delta(c\eta^n))\right)< \infty.
\end{equation}

In respect of (\ref{ecn3}), using that $\{k^{-1}\eee(k)\}_{k\geq 1}$ converges to $\tau$ in a non--increasing way, that $\{W_n\}_{n\geq 0}$ is a
supermartingale, that the process is uniformly supercritical and condition (iii), one has that
$$
\sum_{n=0}^\infty \sigma_n^2\frac{E[\eee(Z_n)]}{r^2_{n+1}}< \sum_{n=0}^\infty \tau\sigma_n^2\frac{E[Z_n]}{(\tau m_n^2)r_{n}}\leq \frac{N}{\tau}\sum_{n=0}^\infty\frac{\sigma_n^2}{m_n^2r_n}< \infty.
$$
The convergence in (\ref{ecn2}) follows from  $\sum_{k=1}^\infty k^{-1}\delta(k)<\infty$ and $\sum_{k=1}^\infty k^{-1}h(k)<\infty$ assumed  in a) and b).

Finally, the $L^2$--convergence of $\{W_n\}_{n\geq 0}$ is proved using Doob's decomposition and  following  similar ideas to those used in Theorem 3 in  \cite{GonzalezMolinaPuertoc}.

\section*{Acknowledgement}
Research supported by the Ministerio de Econom\'ia y Competitividad and the FEDER through the Plan Nacional de Investigaci\'on Cient\'ifica, Desarrollo e Innovaci\'on Tecnol\'olgica, grant MTM2012-31235.

This is a plain preprint version of the article that was published, after peer review, and is subject to Springer Nature’s terms of use, but is not the Version of Record and does not reflect post-acceptance improvements, or any corrections. The Version of Record is available online at: \url{http://dx.doi.org/10.1007/s10986-015-9265-0}.

%\bibliography{x}
%\bibliographystyle{plain}

\newcommand{\nosort}[1]{}

\end{document}